\documentclass[twoside,bezier]{amsart}
\usepackage[left=2cm,right=2cm]{geometry}

\usepackage{amsmath,amsthm,amsfonts,amssymb,amscd,euscript}
\usepackage{setspace}

\usepackage[figuresright]{rotating}
\newtheorem{thm}{Theorem}[section]
\newtheorem{cor}[thm]{Corollary}
\newtheorem{lem}[thm]{Lemma}
\newtheorem{exam}[thm]{Example}
\newtheorem{prop}[thm]{Proposition}
\theoremstyle{definition}\newtheorem{defn}[thm]{Definition}
\theoremstyle{remark}

\numberwithin{equation}{section}
\newcommand{\filebegin}{\begin{document}}
\newcommand{\fileend}{\end{document}}
\makeatother
\def\thefootnote{}
\newcommand{\lo}{\longrightarrow}
\newcommand{\NMM}{\hspace*{2mm}}

\usepackage[left=2cm,right=2cm]{geometry}

\def\n{\noindent}%

\numberwithin{equation}{section}
\def\mapdown#1{\Big\downarrow\rlap
{$\vcenter{\hbox{$\scriptstyle#1$}}$}}
\begin{document}

\vspace*{2cm}
\begin{center}
{\bf\large $C^{*}$-semi-inner product spaces}
 \\[0.5cm]
{Mohammad Janfada, Saeedeh Shamsi Gamchi, Assadollah Niknam\\
Department of Pure Mathematics, Ferdowsi University of Mashhad, Mashhad, P.O. Box 1159-91775, Iran\\

Email: mjanfada@gmail.com\\
Email: saeedeh.shamsi@gmail.com\\
Email: dassamankin@yahoo.co.uk} \\[2mm]

\end{center}%
\vspace*{0.5cm}
\begin{quotation}
\noindent
{\footnotesize
{\sc Abstract.}
In this paper we introduce a generalization of Hilbert $C^*$-modules which are pre-Finsler module namely $C^{*}$-semi-inner product spaces. Some properties and results of such spaces are investigated, specially the orthogonality in these spaces will be considered. We then study bounded linear operators on $C^{*}$-semi-inner product spaces.}
\end{quotation}

\n \textbf{2010 Mathematics subject classification: } 46C50, 46L08, 47A05.%\\

\markboth
{M. Janfada, S. Shamsi Gamchi , A. Niknam}
 {$C^{*}$-semi-inner product spaces}

%%%%%%%%%%%%%%%%%%%%%%%%%%%%%%%%%%%%%%%%%%%%%%%%%%%%%%%%%%%%%%%%%%%%%%%%%%%%%%%%%%%%%%%%%%%%%%%%%%%%%%%%%%%%%%%%%%

%%%%%%%%%%%%%%%%%%%%%%%%%%%%%%%%%%%%%%%%%%%%%%%%%%%%%%%%%%%%%%%%%%%%%%%%%%

\section{ Introduction}
The semi-inner product (s.i.p., in brief) spaces were introduced by Lumer in \cite{Lum1}, he considered vector spaces on which instead of a bilinear form there is defined a form $[x,y]$ which is linear in one component only, strictly positive, and satisfies Cauchy-Schwarz's inequality. Six years after Lumer's work, Giles in \cite{Gil} explored fundamental properties and consequences of semi-inner product spaces. Also, a generalization of semi-inner product spaces was considered by replacing Cauchy-Schwarz's inequality by Holder's inequality in \cite{Nath}. The concept of $\ast$-semi-inner product algebras of type(p) was introduced and some properties of such algebras were studied by Siham Galal El-Sayyad and S. M. Khaleelulla in \cite{Siham}, also, they obtained some interesting results about generalized adjoints of bounded linear operators on semi-inner product spaces of type(p). In the sequel, a version of adjoint theorem for maps on semi-inner product spaces of type(p) is obtained by Endre Pap and Radoje Pavlovic in \cite{Pap2}. The concept of s.i.p. has been proved useful both theoretically and practically. The applications of s.i.p. in the theory of functional analysis was demonstrated, for example, in \cite{Dra1,Dra2,Fau,Koe,Lum2,Put,Tor, Zha}.\\
On the other hand the concept of a Hilbert $C^*$-module which is a generalization of the notion of a
Hilbert space, first made by I. Kaplansky in 1953 (\cite{Kap}). The research on Hilbert $C^*$-modules began in the 70es (W.L. Paschke,
\cite{Pas}; M.A. Rieffel, \cite{Rie}). Since then, this generalization of Hilbert spaces was considered by many mathematicians, for more details about Hilbert $C^*$-modules we refer also to \cite{Lance}. Also Finsler modules over $C^*$-algebras as a generalization of Hilbert $C^*$-modules,  first investigated in \cite{Phil}. For more on Finsler modules, one may see \cite{Amy,Bak}.\\
In this paper we are going to introduce a new generalization of Hilbert $C^*$-modules which are between Hilbert $C^*$-modules and Finsler modules. Furthermore, $C^{*}$-semi-inner product space is a natural generalization of a semi-inner product space arising under replacement of the field
of scalars $\mathbb{C}$ by a $C^*$-algebra.

\section{$C^{*}$-semi-inner product space}
In this section we investigate basic properties of $C^{*}$-semi-inner product spaces.
\begin{defn}
Let $\mathcal{A}$ be a $C^{*}$-algebra and $X$ be a right
$\mathcal{A}$-module. A mapping $[.,.]: X\times X\rightarrow
\mathcal{A}$ is called a $C^{*}$-semi-inner product or
$C^{*}$-s.i.p., in brief, if the following properties are
satisfied:
\\(i) $[x,x]\geq 0$, for all $x\in X$ and $[x,x]=0$
implies $x=0$;\\
(ii) $[x,\alpha y_{1}+\beta y_{2}]= \alpha
[x,y_{1}]+\beta[x,y_{2}]$, for all $x,y_{1},y_{2}\in X$ and
$\alpha, \beta \in \mathbb{C}$;
\\(iii) $[x,ya]=[x,y]a$ and
$[xa,y]=a^*[x,y]$, for all $x,y\in X$ and $a\in \mathcal{A}$;
\\(iv)
$|[y,x]|^{2}\leq \|[y,y]\|[x,x]$.\\
The triple $(X,\mathcal{A}, [.,.])$ is called a $C^*$-semi-inner
product space or we say $X$ is a semi-inner product
$\mathcal{A}$-module.\\
The property (iv) is called the Cauchy-Schwarz inequality.\\
 If $\mathcal{A}$ is a unital $C^*$-algebra, then one may see that $[\lambda x,y]=\overline{\lambda}[x,y]$, for all $x,y\in X$ and $\lambda\in\mathbb{C}$. Indeed, by the property (iii) we have
 $$[\lambda x,y]=[x(\lambda 1),y]=(\lambda 1)^*[x,y]=\overline{\lambda}[x,y].$$
\end{defn}

One can easily see that every Hilbert $C^*$-module is a
$C^*$-semi-inner product space, but the converse is not true in
general. The following is an example of a
$C^*$-semi-inner product space which is not a Hilbert
$C^*$-module. First we recall  that a semi-inner-product (s.i.p.) in the sense of Lumer and Giles  on
a complex vector space $X$ is a complex valued function $[x,y]$ on
$X\times X$ with the following properties:\\
1. $[\lambda y+z,x]= \lambda[y,x]+ [z,x]$ and $[x,\lambda y]=\overline{\lambda}[x,y]$,
for all complex $\lambda$,\\
2. $[x,x]\geq 0$, for all $x\in X$ and $[x,x]=0$ implies $x=0$;\\
3. $|[x,y]|^2\leq[x,x][y,y]$.\\
 A vector space with a s.i.p. is called a semi-inner-product space (s.i.p.
space) in the sense of Lumer-Giles(see \cite{Lum1}). In this case one may prove that $\|x\|:=[x,x]^{\frac{1}{2}}$ define a norm on $X$. Also it is  well-known  that for every Banach space $X$, there exists a semi-inner product whose norm is equal to its original norm.\\
It is trivial that every Banach space is a semi-inner product $\mathbb{C}$-module.
\begin{exam}\label{c}
Let $\Omega$ be a set and let for any $t\in \Omega$, $X_t$
be a semi-inner product space with the semi inner product $[.,.]^{X_t}$. Define
$$[x,y]_{X_t}:= [x,y]^{X_t}, ~~~x, y \in X_t,$$
trivially $[x,\alpha y+z]_{X_t}= \alpha
[x,y]_{X_t}+[x,z]_{X_t}$ and $[\alpha x,y]_{X_t}=\overline{\alpha}[x,y]_{X_t}$. Let $B =\cup_{t}X_t$ be a bundle of
these semi-inner product spaces over $\Omega$. Suppose $\mathcal{A}=Bd(\Omega)$, the set of all bounded complex-valued functions on $\Omega$,  and $X$ is the
set of all maps $f :
\Omega\to B$ such that $f(t) \in X_t$, for any $t\in \Omega$, with $\sup_{t\in \Omega}\|f(t)\|<\infty$. One can easily see that $X$ is naturally a
$Bd(\Omega)$-module. Furthermore it has a $Bd(\Omega)$-valued semi-inner
product defined by
\[[f, g](t) = [f(t), g(t)]_{X_t},\]
for $t\in \Omega$, hence, it is a $C^*$-semi-inner product space. One
can easily verify that the properties of $C^*$-semi-inner product
are valid.
\end{exam}
Suppose $(\mathcal{A}_{i},\|.\|_{i})$'s, $1\leq i\leq n$, are $C^*$-algebras, then $\bigoplus_{i=1}^{n}\mathcal{A}_{i}$ with its point-wise operations is a $C^*$-algebra  Moreover, $\|(a_1,...,a_n)\|=\max_{1\leq i\leq n} \|a_i\|$ is a $C^*$-norm on $\bigoplus_{i=1}^{n}\mathcal{A}_{i}$. Note that $(a_1,...,a_n)\in(\bigoplus_{i=1}^{n}\mathcal{A}_{i})_+ $ if and only if $a_i \in (\mathcal{A}_i)_+$. Now we may construct the following example.\\
\begin{exam} Let $(X_i,[ . , .  ]_i)$ be a semi-inner product $\mathcal{A}_i$-module, $1\leq i\leq n$. If for $(a_1,...,a_n)\in \mathcal{A}$ and $(x_1,...,x_n)\in \bigoplus_{i=1}^{n}X_{i}$, we define $(x_1,...,x_n)(a_1,...,a_n)=(x_1 a_1,...,x_n a_n)$ and the $C^*$-s.i.p. is defined as follows
 $$[(x_1,...,x_n),(y_1,...,y_n)]=([x_1,y_1]_1,...,[x_n,y_n]_n)$$
  then the direct sum  $\bigoplus_{i=1}^{n}X_{i}$ is a semi-inner product $\mathcal{A}$-module, where $\mathcal{A}=\bigoplus_{i=1}^{n}\mathcal{A}_{i}$.
  \end{exam}

Let $(X,\mathcal{A}, [.,.])$ be a $C^*$-semi-inner product space.
For any $x\in X$, put $|||x|||:=\|[x,x]\|^\frac{1}{2}$. The
following proposition shows that $(X, |||. |||)$ is a normed
$\mathcal{A}$-module.
\begin{prop}
Let $X$ be a right $\mathcal{A}$-module and $[.,.]$ be a
$C^{*}$-s.i.p. on $X$. Then the mapping $x\rightarrow \|[x,
x]\|^{\frac{1}{2}}$ is a norm on $X$. Moreover, for each $x\in X$
and $a\in \mathcal{A}$ we have
$|||xa|||\leq |||x||| \ \|a\|$.
\end{prop}
\begin{proof} Clearly $|||x|||=\|[x,x]\|^\frac{1}{2}\geq
0$ and $|||x|||=0$ implies that $x=0$.\\
Also for each $x\in X$, $\lambda\in \mathbb{C}$, by the Cauchy-Schwarz inequality,
\begin{align*}
|||\lambda x|||^{2}= \|[\lambda x,\lambda
x]\|&=|\lambda|\|[\lambda x, x]\|\\&\leq|\lambda| \ \| \ |[\lambda
x, x]| \ \|\\&\leq|\lambda| \ |||\lambda x||| \ |||x|||.
\end{align*}
Hence, $|||\lambda x|||\leq |\lambda| \ |||x|||$. On the other
hand, we have $|||x|||=|||\frac{1}{\lambda}.\lambda x|||\leq
\frac{1}{|\lambda|}|||\lambda x|||$, therefore, $|||\lambda
x|||=|\lambda| \ |||x|||$.\\
Finally for each $x,y\in X$,
\begin{align*}
|||x+y|||^{2} = \|[x+y,x+y]\| &\leq \|[x+y,x]\|+\|[x+y,y]\|
\\&\leq \| \ |[x+y,x]| \ \|+\| \ |[x+y,y]| \ \| \\&\leq |||x+y|||
\ |||x|||+ |||x+y||| \ |||y||| \\&\leq |||x+y|||(|||x|||+|||y|||).
\end{align*}
Therefore, $|||x+y|||\leq |||x|||+|||y|||$.\\
Also we have
\begin{align*}
|||xa|||^{2}=\|[xa,xa]\|&=\|[xa,x]a\|\\&\leq \|[xa,x]\| \
\|a\|\\&\leq |||xa||| \ |||x||| \ \|a\|,
\end{align*}
hence, $|||xa|||\leq |||x||| \ \|a\|$.\\
\end{proof}
As another result for this norm one can see that for each  $x\in X$,   $|||x[x,x]|||=|||x|||^{3}$. Indeed,
\begin{eqnarray*}
|||x[x,x]|||^{2}&=&\|[x[x,x],x[x,x]]\| \\
&=&\| [x,x]^{3}\| \\
&=&\|[x,x]\|^3.
\end{eqnarray*}
The last equality follows from the fact that in any $C^*$-algebra, we have $\|a^3\|=\|a\|^3$, for any self-adjoint element $a\in
\mathcal{A}$.
\begin{prop}
Let $\mathcal{A}$ and $\mathcal{B}$ be two $C^*$-algebras and
$\psi: \mathcal{A} \to \mathcal{B} $ be an $\ast$-isomorphism. If
$(X,[.,.]_{\mathcal{A}})$ is a $C^*$-semi-inner product
$\mathcal{A}$-module, then $X$ can be represented as a right
$\mathcal{B}$-module with the  module action $x\psi(a)=xa$ and is a $C^*$-semi-inner product
$\mathcal{B}$-module with the  $C^*$-semi-inner product defined by
 $$[.,.]_ \mathcal{B} =\psi ([.,.]_\mathcal{A}). $$
\end{prop}
\begin{proof} It is clear that $X$ is a right $\mathcal{B}$-module
with the mentioned module product.
It is easy to verify that the properties (i) to (iii) of definition
of $C^*$-semi-inner product holds for $[.,.]_\mathcal{B}$. Now, we prove the property (iv) for
$[.,.]_\mathcal{B}$. Since $\psi: \mathcal{A} \to \mathcal{B} $ is
an $\ast$-isomorphism, so it is isometric and $\psi
(\mathcal{A}_+)\subseteq \mathcal {B}_+$. Thus we have
\begin{align*}
|[x,y]_\mathcal{B}|^2
=|\psi([x,y]_\mathcal{A})|^2&=\psi([x,y]_\mathcal{A})^*
\psi([x,y]_\mathcal{A})\\&=\psi([x,y]_\mathcal{A}^*[x,y]_\mathcal{A})\\&=\psi(|[x,y]_\mathcal{A}|^2)\\&\leq
\|[x,x]_\mathcal{A}\| \
\psi([y,y]_\mathcal{A})\\&=\|\psi([x,x]_\mathcal{A})\| \
\psi([y,y]_\mathcal{A})\\&=\|[x,x]_\mathcal{B}\| \
[y,y]_\mathcal{B}.
\end{align*}
\end{proof}
We will establish a converse statement to the above proposition. Consider that a semi-inner product $\mathcal{A}$-module $X$ is said to be full if the linear span of $\{[x,x]: x\in X\}$, denoted by $[X,X]$, is dense in $\mathcal{A}$.
\begin{thm}
Let $X$ be both a full complete semi-inner product $\mathcal{A}$-module and a full complete semi-inner product $\mathcal{B}$-module such that $\|[x,x]_{\mathcal{A}}\|=\|[x,x]_{\mathcal{B}}\|$ for each $x\in X$, and let $\psi: \mathcal{A} \to \mathcal{B} $ be a map such that $xa=x\psi(a)$ and $\psi ([x,x]_\mathcal{A})=[x,x]_ \mathcal{B} $ where $x\in X, \ a\in\mathcal{A}$. Then $\psi$ is an $\ast$-isomorphism of $C^*$-algebras.
\end{thm}
\begin{proof}
The proof is similar to theorem $2.1$\cite{Amy}.
\end{proof}
We recall that if $\mathcal{A}$ is a C*-algebra, and $\mathcal{A}_+$ is the set
of positive elements of $\mathcal{A}$, then  a pre-Finsler
$\mathcal{A}$-module is a right $\mathcal{A}$-module $E$ which is
equipped with a map $\rho : E \to \mathcal{A}_+$ such that\\
(1) the map $\| . \|_E : x \mapsto \|\rho(x)\|$ is a norm on $E$; and\\
(2) $\rho(xa)^2 = a^*\rho(x)^2 a$, for all $a\in \mathcal{A}$ and $x\in E$.\\
If $(E,\|.\|_E)$ is complete then $E$ is called a Finsler
$\mathcal{A}$-module. This definition is a modification of one
introduced by N.C. Phillips and N.Weaver \cite{Phil}. Indeed it is
routine by using an interesting theorem of C. Akemann
[\cite{Phil}, Theorem 4] to show that the norm completion of a
pre-Finsler $\mathcal{A}$-module is a Finsler
$\mathcal{A}$-module. Now it is trivial to see that every
$C^*$-semi-inner product space $(X,\mathcal{A}, [.,.])$ is a
pre-Finsler module with the function $\rho :X\to \mathcal{A}_+$
defined by $\rho(x)=[x,x]^{\frac{1}{2}}$. Thus every complete
$C^*$-semi-inner product space enjoys all the properties of a
Finsler module.
\begin{prop}\cite{Phil} Let $\mathcal{A}=C_0(X)$ and let $E$ be a
Finsler $\mathcal{A}$-module. Then $\rho$ satisfies  \[\rho(x + y)
\leq \rho(x)+\rho(y)\] for all $x,y\in E$
\end{prop}
Replacing the real numbers, as the codomain of a norm, by an ordered Banach space
we obtain a generalization of normed space. Such a generalized
space, called a cone normed space, was introduced by Rzepecki \cite{R}.
\begin{cor} Let $(X,[.,.])$ be a semi-inner $C(X)$-module, then
$\|.\|_{c} : X \to C(X)$ defined by
$\|x\|_{c}=[x,x]^\frac{1}{2}$ is a cone norm on $X$.
\end{cor}

\section{Orthogonality in $C^{*}$-semi-inner product spaces}
In this section we study the relations between Birkhoff-James orthogonality and the orthogonality in a $C^{*}$-semi-inner product spaces.\\
In a normed space $X$ (over $\mathbb{K}\in\{\mathbb{R},\mathbb{C}\}$), the Birkhoff-James orthogonality
(cf.\cite{B,J}) is defined as follows
$$x\perp_{B} y \ \Leftrightarrow \
\forall\alpha\in \mathbb{K}; \  \|x+\alpha y\|\geq \|x\|.$$

\begin{thm}
Let $X$ be a right $\mathcal{A}$-module and $[.,.]$ be a
$C^{*}$-s.i.p. on $X$. If $x,y\in X$ and $[x,y]=0$ then $x\perp_{B} y$.
\end{thm}
\begin{proof} Let $[x,y]=0$.  If $x=0$ then by the definition of
Birkhoff-James orthogonality it is obvious that $x\perp_{B} y$.
Now if $x\neq 0$, then for all $\alpha\in \mathbb{K}$,
\begin{align*}
|||x|||^2- |\alpha| \ \|[x,y]\|&\leq\|[x,x+\alpha y]\|\\&\leq
|||x||| \ |||x+\alpha y|||.
\end{align*}
Hence, $$-|\alpha| \ \|[x,y]\|\leq |||x|||(|||x+\alpha
y|||-|||x|||).$$
But $x\neq 0$ and $[x,y]=0$, so by the above
inequality we conclude that $|||x+\alpha y|||\geq |||x|||$, which
shows that $x\perp_{B} y$.
\end{proof}
In the sequel we try to find a sufficient condition for $x,y$ to be orthogonal in the $C^*$-semi-inner product. For; we need some preliminaries. we remind that in a $C^*$-algebra $\mathcal{A}$ and for any $a\in \mathcal{A}$ there exist self-adjoint elements $h,k\in \mathcal{A}$ such that $a=h+ik$. We apply $Re(a)$ for $h$.
\begin{defn}
A $C^{*}$-s.i.p. $[.,.]$ on right $\mathcal{A}$-module $X$ is said
to be continuous if for every $x,y\in X$ one has the
equality
$$\lim_{t\rightarrow 0} Re[x+ty,y] = Re [x,y],$$
where $t\in\mathbb{R}$.
\end{defn}
\begin{exam}
In Example $\ref{c}$, $\Omega=\{1,2,...,n\}$ and $X$ be the  semi inner product $Bd(\Omega)$-module defined in Example $\ref{c}$. If $X_t$ is a continuous s.i.p. space (see \cite{Gil}), for all $t\in \Omega$, then $X$ is a continuous $C^*$-s.i.p space. Indeed, it is clear that
$$\sup_{t\in \Omega}\|Re[f(t)+\alpha g(t),g(t)]_{X_t}-Re[f(t),g(t)]_{X_t}\|$$ tends to $0$, when $\alpha \rightarrow 0$.
\end{exam}
\begin{thm} Let $X$ be a right $\mathcal{A}$-module and let $[.,.]$ be a continuous $C^{*}$-s.i.p. on
$X$ such that $[x,y]\in \mathcal{A}_{sa}$ for each $x,y\in X$. If for $x,y\in X$ and any $t\in
\mathbb{R}$,
$$[x+ty,x+ty]\geq
[x,x]^\frac{1}{2} \ |||x+ty|||$$
then $[x,y]=0$.
\end{thm}
\begin{proof} It is clear that for each $a\in \mathcal{A}_{sa}$, we have $a\leq |a|$. Now assume that
$$[x+ty,x+ty]\geq [x,x]^\frac{1}{2} \ |||x+ty|||$$for all $x,y\in
X$ and $t\in \mathbb{R}$. By Cauchy-Schwarz inequality (iv) and the
fact that $[x,y]\in \mathcal{A}_{sa}$ for each $x,y\in X$, we get;
\begin{align*}
[x+ty,x+ty]&\geq [x,x]^\frac{1}{2} \ |||x+ty|||\\&\geq
|[x+ty,x]|\\&\geq [x+ty,x]
\end{align*}
so, we have: $t[x+ty,y]\geq 0$ for each $t\in\mathbb{R}$. Thus,
for $t\geq 0$ we have $[x+ty,y]\geq 0$ and for $t\leq 0$ we have
$[x+ty,y]\leq 0$. Now, since $[.,.]$ is a continuous
$C^{*}$-s.i.p. and $\mathcal{A}_{+}$ is a closed subset of
$\mathcal{A}$, so we have
\begin{align*}
 0\geq [x,y]&=\lim_{t\rightarrow 0^{-}}[x+ty,y] \\&= \lim_{t\rightarrow
0^{+}} [x+ty,y] = [x,y]\geq 0,
\end{align*}
thus, $[x,y]=0$.
\end{proof}

\section{Bounded Linear operators on $C^{*}$-semi-inner product spaces}

\begin{thm}
Let $X$ be a semi inner product $\mathcal{A}$-module. Then for
every $y\in X$ the mapping $f_{y}: X\rightarrow \mathcal{A}$
defining by $f_{y}(x)=[y,x]$ is a $\mathcal{A}$-linear continuous
operator endowed with the norm generated by
$C^{*}$-s.i.p. Moreover, $\|f_{y}\|=|||y|||$.
\end{thm}
\begin{proof} The fact that $f_{y}$ is a $\mathcal{A}$-linear
operator follows by (ii) and (iii) of definition $1.1$. Now, using
Schwartz's inequality (iv) we get;
$$\|f_{y}(x)\|=\|[y,x]\|\leq |||y||| \ |||x|||$$
which implies that $f_{y}$ is bounded and $$\|f_{y}\|\leq
|||y|||$$On the other hand, we have;$$||f_{y}||\geq
\|f_{y}(\frac{y}{|||y|||})\|=|||y|||$$and then
$\|f_{y}\|=|||y|||$.
\end{proof}

\begin{cor}
If  $X$ is a right $\mathcal{A}$-module and $[.,.]$ a
$C^{*}$-s.i.p. on X, then for all $x\in X$ we have; $$|||x||| =
\sup \{ \|[x,y]\| : \ |||y|||\leq 1\}.$$
\end{cor}

\begin{lem}
\cite{Jo,P} Let $\mathcal{A}$ be a unital $C^*$-algebra let
$r:\mathcal{A}\rightarrow \mathcal{A}$ be a linear map such that
for some constant $K\geq0$ the inequality $r(a)^*r(a)\leq Ka^*a$
is fulfilled for all $a\in \mathcal{A}$. Then $r(a)=r(1)a$ for all
$a\in \mathcal{A}$.
\end{lem}

\begin{thm}
 Let $X$ and $Y$ be semi inner product $\mathcal{A}$-modules, $T:X\rightarrow Y$ be a linear
 map. Then the following conditions are equivalent:

(i) the operator $T$ is bounded and $\mathcal{A}$-linear, i,e.
$T(xa)= Tx.a$ for all $x\in X$, $a\in \mathcal{A}$;

(ii) there exists a constant $K\geq 0$ such that for all $x\in X$
the operator inequality $[Tx,Tx]\leq K[x,x]$ holds.
\end{thm}

\begin{proof}
To obtain the second statement from the first one, assume that $T(xa)= Tx.a$ and $\|T\|\leq 1$. If $C^*$-algebra
$\mathcal{A}$ does not contain a unit, then we consider
modules $X$ and $Y$ as modules over $C^*$-algebra $\mathcal{A}_1$,
obtained from $\mathcal{A}$ by unitization. For $x\in X$ and $n\in
\mathbb{N}$, put

$$a_n=([x,x]+\frac{1}{n})^{-\frac{1}{2}}, \  \  \  \  x_n=x a_n $$
Then $[x_n,x_n]=a_n^*[x,x]a_n=[x,x]([x,x]+\frac{1}{n})^{-1}\leq1$,
therefore, $\|x_n\|\leq 1$, hence $\|Tx_n\|\leq 1$. Then for all
$n\in \mathbb{N}$ the operator inequality $[Tx_n,Tx_n]\leq1$ is
valid. But
$$[Tx,Tx]=a^{-1}_n[Tx_n,Tx_n]a^{-1}_n\leq
a^{-2}_n=[x,x]+\frac{1}{n}.$$Passing in the above inequality to
the limit $n\rightarrow \infty$, we obtain $[Tx,Tx]\leq [x,x]$. To
derive the first statement from the second one we assume that for
all $x\in X$ the inequality $[Tx,Tx]\leq [x,x]$ is fulfilled and it
obviously implies that the operator $T$ is bounded,
$\|T\|\leq 1$. Let $x\in X$, $y\in Y$. Let us define a map
$r:\mathcal{A}_1\rightarrow\mathcal{A}_1$ by the equality

$$r(a)=[y,T(xa)].$$
Then

$$r(a)^*r(a)=|[y,T(xa)]|^2 \leq |||y|||^2[T(xa),T(xa)]\leq
|||y|||^2[xa,xa]=|||y|||^2a^*[x,x]a\leq
|||y|||^2|||x|||^2a^*a.$$Therefore, by the above lemma we have
$r(a)=r(1)a$, i.e.$$[y,T(xa)]=[y,Tx]a =[y,Tx.a]$$for all
$a\in\mathcal{A}$ and all $y\in Y$. Hence, the proof is complete.

\end{proof}
\begin{cor}
Let $X$ and $Y$ be semi inner product $\mathcal{A}$-modules,
$T:X\rightarrow Y$ be a bounded $\mathcal{A}$-linear
 map. Then $$\|T\|=\inf \{ K^{\frac{1}{2}}: [Tx,Tx]\leq K[x,x] \}.$$
\end{cor}

%========================================================
%===================================================
%\section {\sc }

%\section {\sc }

%=======================================================

\end{document}